\documentclass[12pt]{iopart}
\usepackage{amsmath}
\usepackage{iopams}
\usepackage{graphicx}
\usepackage{float}
\restylefloat{figure}

\def\RR{\mathbb R}
\newcommand{\qed}{\hbox to 0pt{}\hfill$\rlap{$\sqcap$}\sqcup$\vspace{2mm}}
\newtheorem{theorem}{Theorem}[section]

\numberwithin{equation}{section}

\newtheorem{guess}[theorem]{Theorem}
\newtheorem{remark}[theorem]{Remark}
\newtheorem{corol}[theorem]{Corollary}
\newtheorem{example}[theorem]{Example}
\newtheorem{definition}[theorem]{Definition}

\begin{document}

\title[Stability of cooperative systems with distributed delays]{On stability of 
cooperative and hereditary systems with a distributed delay}

\author{Leonid Berezansky$^1$ and Elena Braverman$^2$}
\address{$^1$ Dept. of Math, Ben-Gurion University of Negev, Beer-Sheva 84105, 
Israel}
\ead{brznsky@cs.bgu.ac.il}
\address{$^2$ Dept. of Math \& Stats, University of Calgary,
2500 University Dr. NW, Calgary, AB, Canada T2N 1N4}
\ead{maelena@math.ucalgary.ca}


\begin{abstract}
We consider a system $\displaystyle \frac{dx}{dt}=r_1(t) G_1(x) \left[ \int_{h_1(t)}^t f_1(y(s))~d_s R_1 (t,s) - x(t) 
\right], \\ \frac{dy}{dt}=r_2(t) G_2(y)  \left[ \int_{h_2(t)}^t f_2(x(s))~d_s R_2 (t,s) - y(t)\right]$ 
with increasing functions $f_1$ and $f_2$, which has at most one positive equilibrium.  
Here the values of the functions $r_i,G_i,f_i$ are positive for positive arguments,
the delays in the cooperative term can be 
distributed and unbounded, both systems with concentrated delays and integro-differential systems are a particular 
case of the considered system.
Analyzing the relation of the functions $f_1$ and $f_2$, we obtain several possible scenarios of the global
behaviour. They include the cases when all nontrivial positive solutions tend to the same attractor which can be
the positive equilibrium, the origin or infinity. 
Another possibility is the dependency of asymptotics on the initial conditions: either
solutions with large enough 
initial values tend to the equilibrium, while others tend to zero, or solutions with small enough
initial values tend to the equilibrium, while others infinitely grow.
In some sense solutions of the equation are intrinsically non-oscillatory:
if both initial functions are less/greater than the equilibrium value,
so is the solution for any positive time value. The paper continues the study of equations with monotone 
production functions initiated in [Nonlinearity, 2013, 2833-2849].
\end{abstract}



\maketitle

\section{Introduction}

The system of autonomous differential equations with constant delays in the production term
\begin{equation}
\label{intro1}
\begin{array}{ll}
\displaystyle \frac{dx}{dt} & = R_1(y(t-\tau_1)) - a_1x(t)  \vspace{2mm} \\
\displaystyle \frac{dy}{dt} & = R_2(x(t-\tau_2)) - a_2 y(t)
\end{array}
\end{equation}
was considered in \cite{Vargas2}, where $R_i:(0,\infty)\to (0,\infty)$ are monotone increasing functions. 
It can describe a couple of populations, where the growth of each population is stimulated by 
the size of the other population and is suppressed by its own growth. 
Systems of differential equations describing different types of species, where the rate of change for
each of them is positively influenced by all other populations but itself, are usually called {\em cooperative}.
This is in contrast, for example, to competitive systems, where this influence is negative, and predator-prey systems,
with different types of influences. These systems can correspond to  the cooperative types of species, or to the patch 
environment, 
the growth in each patch is suppressed by overpopulation in itself while stimulated by high density in adjacent patched,
due, for example, to possible immigration.
Another situation is {\em hereditary systems} where each variable describes a different developmental stage of the 
same species (e.g. eggs, larvae, juveniles, adult species capable of reproduction). In the case of system \eqref{intro1}, 
$x$ and $y$ can be juvenile and adult counts, respectively. There is a competition within each group, as well as natural
mortality, and the mortality per capita rate is assumed to be population-independent.
All the growth of juveniles is due to reproduction of adults, while maturation of juveniles contributes to adult 
numbers. There are delays in both recruitment processes (maturation delay for juveniles and reproduction time for adults).
In line with the above description, model (\ref{intro1}) includes delay in the reproduction term only, and the mortality is 
assumed to be proportional to the current population density.

In the present paper, we consider systems of two equations where the growth of each of two variables is stimulated by
high numbers in the other (due to cooperation, or inheriting part of it, or influx of offspring of the other population),
and call them cooperative or hereditary systems. The delays of a positive impact can describe the time required to 
translate nutritional benefits into body mass for the cooperation type. For hereditary
systems, we have maturation and reproduction delays.

System (\ref{intro1}) includes the two-neuron bidirectional associative memory (BAM) model \cite{Kosko} 
\begin{equation}  
\label{intro1a}
x'(t)=-x(t)+af(y(t))+I,~~y'=-y(t)+bg(x(t))+J.
\end{equation}
A simplified version of the delay system considered in \cite{Dong}
\begin{equation}
\label{intro6}
\begin{array}{ll}
\displaystyle \frac{dx}{dt} & = c_1 \tanh(y(t-\tau_1)) - \mu_1 x(t)   \vspace{2mm} \\
\displaystyle \frac{dy}{dt} & = c_2 \tanh(x(t-\tau_2)) - \mu_2 y(t)
\end{array}
\end{equation}
is also a particular case of (\ref{intro1}).

Another autonomous model
\begin{equation}
\label{intro2}  
\begin{array}{ll}
\displaystyle \frac{dx}{dt} & \displaystyle  = G_1(x(t)) \left[ R_1(y(t-\tau_1)) - a_1x(t) \right]  \vspace{2mm} \\
\displaystyle \frac{dy}{dt} & \displaystyle = G_2(y(t)) \left[ R_2(x(t-\tau_2)) - a_2 y(t) \right]
\end{array}
\end{equation}  
includes a system of logistic equations with the delay in the production term; equations of this type were described in 
\cite{AWW}. Some particular non-delay systems of type (\ref{intro2}) were studied in \cite{Vargas1}.
For example, the Lotka-Volterra cooperative model considered in \cite{LuLu,Saito,add1}, if the delayed mortality 
terms are omitted, has the form
\begin{equation}
\label{intro3}
\begin{array}{ll}
\displaystyle \frac{dx}{dt} & \displaystyle = x(t) \left[ r_1 - a_1 x(t) + b_{1} y(t-\tau_1) \right]  \vspace{2mm} \\
\displaystyle \frac{dy}{dt} & \displaystyle = y(t) \left[ r_2 - a_2 y(t) + b_{2} x(t-\tau_2) \right]~.  
\end{array}
\end{equation}
Evidently \eqref{intro3} is a particular case of (\ref{intro2}), and all the results of \cite{Vargas2} are applicable 
to (\ref{intro3}). 

The Hopfield neural network \cite{Hopfield} 
\begin{equation} 
\label{intro2a}
x_i^{\prime}(t) = -b_i(x_i(t))+ \sum_{j=1}^n c_{ij} f_j(x_j(t))+I_i, ~~t\geq 0,~i=1, \dots, n,
\end{equation}
with $n=2$ and $c_{ii}=0$, can be rewritten as (\ref{intro2}) with $\tau_1=\tau_2=0$, arbitrary $a_i>0$ and
$G_i(u) = b_i(u)/(a_i u)$, $R_i(u)= c_{ij}f_j(u)/G_i(u)$, $j \neq i$.

The purpose of the present paper is to explore global asymptotic stability of 
cooperative systems with a distributed delay, which include (\ref{intro1}) and (\ref{intro2}) as special cases;
in addition to being distributed, the delay can change with time.
Distributed delays describe a feasible fact that any interval for delay value has some probability,
such models include equations with concentrated (either constant or variable) delays. Stability
of equations and systems with distributed delays  attracted 
recently much attention, see, for example, \cite{Automatica2012,Nonlin2013,JMAA2014,BrZhuk,
Oliveira,FariaOliveira, FariaTrof, Kyrychko, Beretta,Liu,Ncube,Niu,Ngoc,Park,Solomon,Yuan} 
for some recent results and their applications, also see 
references therein. The summary of the results obtained by the beginning of 1990ies can be found in \cite{Kuang}.
The methods applied to establish absolute convergence of the system either to the origin, or to the unique positive 
equilibrium, or to infinity, goes back to \cite{Damir,BrZhuk} and was applied in \cite{Nonlin2013,JMAA2014}.
In contrast to our earlier papers \cite{Damir,BrZhuk,Nonlin2013,JMAA2014}, in the present paper we consider a system, not a 
single equations.
Compared to all other previous work, the main differences are outlined below.
\begin{itemize}
\item
We consider distributed delays of the most general type; as particular cases, they include
systems with variable concentrated delays, integral terms (in most papers, distributed delay is associated
with these integral terms), their combination, and some other models (for example, Cantor function as a distribution).
Moreover, argument deviations can be Lebesgue measurable functions, they are not required to be continuous.
Thus the methods developed for continuous delays are not applicable in this setting.
\item
The delay distributions can be non-autonomous. If we describe these distributions as a probability that
a delay takes a greater than a given value, this corresponds to time-dependent delay. In applications,
this allows to consider, for example, seasonal changes in delay distributions. To some extent, we explore
the most general system with a unique positive equilibrium, and justify global stability of this equilibrium,
once delays are involved in those terms only which describe cross-influences. This is a generalization
of the result in \cite{Vargas2} for a system of two autonomous equations with constant concentrated delays. 
To some extent, we have answered the question when delays do not have any destabilizing effect on a non-autonomous
system of two equations.
\item
On the other hand, many of the previous papers on distributed delay describe much more complicated dynamics than
absolute global stability established in the present paper. For example, delay dependence of stability properties  
was studied in \cite{BrZhuk}, while possible multistability considered in \cite{JMAA2014}.
However, the study of systems which can be destabilized by large enough delay are not in the framework of the present paper.
Here we restrict ourselves to monotone increasing production functions, which can be treated as positive feedback
in the delayed term.
\end{itemize}

The paper is organized as follows.
Section~\ref{permanence} contains existence, positivity and permanence
results for models with a distributed delay. Section~\ref{stability} presents the global
stability theorem which is the main result of the present paper.
Finally, Section~\ref{summary} considers applications and involves some discussion.

\section{Positivity and Solution Bounds}
\label{permanence}

In the present paper we consider the system with a distributed delay 
\begin{equation}
\label{1a}
\begin{array}{ll}
\displaystyle \frac{dx}{dt} & \displaystyle  = r_1(t) \left[ \int_{h_1(t)}^t f_1(y(s))~d_s R_1 (t,s) - x(t) \right]  
\vspace{2mm} \\
\displaystyle \frac{dy}{dt} & \displaystyle  = r_2(t) \left[ \int_{h_2(t)}^t f_2(x(s))~d_s R_2 (t,s) - y(t) \right]
\end{array}
\end{equation}
with the initial conditions 
\begin{equation}
\label{2star}
x(t)=\varphi(t),~t \leq 0, ~~y(t)=\psi(t),~t \leq 0,
\end{equation}
where $\varphi(t)$ and $\psi(t)$ are initial functions.

\begin{definition}
The pair of functions $(x(t),y(t))$ is {\bf a solution  of system} \eqref{1a},\eqref{2star} if
it satisfies \eqref{1a} for almost all $t \geq 0$ and \eqref{2star} for $t \leq 0$.
\end{definition}

System \eqref{1a} will be investigated under some of the following assumptions:
\begin{description}
\item{{\bf (a1)}} 
$f_i:\RR^+ \to \RR^+=[0,\infty)$, $i=1,2$ are continuous functions,  
$f_i$ are strictly monotone increasing on $\RR^+$ ($f_i(x)>f_i(y)$ for $x>y\geq 0$)  and $f_i(x)>0$ for $x>0$, $i=1,2$;
\item{{\bf (a2)}} 
The equation $f_1^{-1}(x)=f_2(x)$ has exactly one positive solution $K>0$, where $f_2(x)>f_1^{-1}(x)$
for $f_1(0)<x<K$ and $f_2(x)<f_1^{-1}(x)$ for $x>K$;
\item{{\bf (a3)}} 
$h_i:\RR^+ \rightarrow \RR$, $i=1,2$ are Lebesgue measurable
functions, $ h_i(t)\leq t$,
$\lim\limits_{t\rightarrow \infty}h_i(t)=\infty$, $i=1,2$;
\item{{\bf (a4)}} 
$R_i(t, \cdot)$, $i=1,2$ are left continuous non-decreasing functions
for any $t$, $R_i(\cdot,s)$ are locally integrable for
any $s$, $R_i(t,s)=0$, $s \leq h_i(t)$, $R_i(t,t^+)=1$, $r_i(t)$ are Lebesgue measurable essentially bounded on $\RR^+$
functions, $r_i(t) \geq 0$, $i=1,2$;
here $u(t^+)$ is the right-side limit of function $u$ at point $t$.
\item{{\bf (a5)}}
${\displaystyle 
\int_0^{\infty} r_i(s)~ds = \infty}, ~ i=1,2$;
\item{{\bf (a6)}} 
$\varphi: (-\infty,0] \to \RR$ and $\psi: (-\infty,0] \to \RR$ are continuous bounded functions,
$\varphi(t) \geq 0$, $\psi(t) \geq 0$, $t <0$, $\varphi(0)>0$, $\psi(0)>0$.
\end{description}

Condition (a2) implies that system (\ref{1a}) has one and only one positive equilibrium 
which is $(x(t),y(t))=(K,f_2(K))$.

As particular cases, system (\ref{1a}) includes the model with variable delays
\begin{equation}
\label{1avar}
\begin{array}{ll}
\displaystyle \frac{dx}{dt} & \displaystyle  = r_1(t) \left[ f_1(y(h_1(t))) - x(t) \right]  
\vspace{2mm} \\
\displaystyle \frac{dy}{dt} & \displaystyle  = r_2(t) \left[ f_2(x(h_2(t))) - y(t) \right]
\end{array}
\end{equation}
where instead of (a4) we assume
\begin{description}
\item{{\bf (b4)}}
$r_i(t)$ are Lebesgue measurable essentially bounded on $\RR^+$
functions, $r_i(t) \geq 0$, $i=1,2$,
\end{description}
and the integro-differential system 
\begin{equation}
\label{1aintegro}
\begin{array}{ll}
\displaystyle \frac{dx}{dt} & \displaystyle  = r_1(t) \left[ \int_{h_1(t)}^t K_1 (t,s)f_1(y(s))~ds  - x(t) \right]
\vspace{2mm} \\  
\displaystyle \frac{dy}{dt} & \displaystyle  = r_2(t) \left[ \int_{h_2(t)}^t K_2 (t,s)f_2(x(s))~ds  - y(t) \right]
\end{array}
\end{equation}
where instead of (a4) we consider the condition
\begin{description}
\item{{\bf (c4)}}   
$K_i(t, s):\RR^+ \times \RR^+ \to \RR^+$, $i=1,2$ are locally integrable functions
in both $t$ and $s$ satisfying $\displaystyle \int_{h_i(t)}^t K_i (t,s)~ds \equiv 1$, 
$r_i(t)$ are Lebesgue measurable essentially bounded on $\RR^+$
functions, $r_i(t) \geq 0$, $i=1,2$.
\end{description}

\begin{definition}
The solution $(x(t),y(t))$ of \eqref{1a},\eqref{2star} is {\bf
permanent} if there exist $a$, $b$ and $A$, $B$, $A\geq a>0$, $B\geq b>0$,  such that
$$a \leq x(t) \leq A, ~~b\leq y(t) \leq B, ~t \geq 0.$$
\end{definition}

Theorem~\ref{theorem1} presents sufficient conditions when there exists a positive solution of \eqref{1a},\eqref{2star} 
on $[0,\infty)$.

\begin{figure}[ht]
\centering
\includegraphics[scale=0.85]{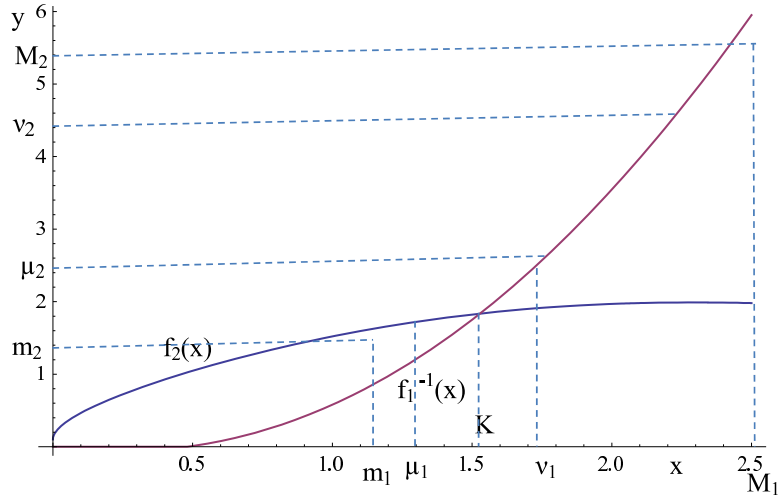} 
\caption{Illustration of the solution bounds}
\label{figure1}   
\end{figure}

\begin{guess}
\label{theorem1}
Suppose (a1),(a3)-(a4),(a6) hold. 

1) A solution of \eqref{1a},\eqref{2star} is positive in its maximal interval of existence $[0,d)$.

2) If in addition 

{\bf (a2$^{\ast}$)} there exists $K>0$ such that $f_2(x)<f_1^{-1}(x)$ for $x>K$

then there exists a positive solution of \eqref{1a},\eqref{2star} for $t \in [0,\infty)$. We will call it a global 
solution.

3) If (a1)-(a4), (a6) hold, then the global solution of \eqref{1a},\eqref{2star} is permanent.
\end{guess}
{\bf Proof.} The proof is illustrated by Fig.~\ref{figure1}.

1) The existence of a local solution which is positive on $[0, \varepsilon)$ is justified 
in the same way as in \cite{Nonlin2013,JMAA2014}, using the result of \cite[Theorem 4.5, p. 95]{Cordun}.

This solution is either global or there exists $t_2$ such that either
\begin{equation}
\label{inf1}
\liminf_{t \to t_2^-} x(t) = -\infty
\end{equation}
or
\begin{equation}
\label{inf2}
\limsup_{t \to t_2^-} x(t) =\infty \,
\end{equation}
or either (\ref{inf1}) or (\ref{inf2}) is satisfied with $y(t)$ instead of $x(t)$.

The initial value is positive, so as long as $x(t)>0$, $y(t)>0$, each component
of the solution $(x(t),y(t))$ is not less than the solution of the initial value problem for the system
of ordinary differential equations
\begin{equation}
\label{th0_1}
x^{\prime}(t)+r_1(t)x(t)=0, ~~ y^{\prime}(t)+r_2(t)y(t)=0, ~~x(t_0)=x_0>0, ~~y(t_0)=y_0>0,
\end{equation}
and this solution is positive for any $t\geq 0$.
Let us assume that either $x(t)$ or $y(t)$ becomes negative and let
$t_1$ be the smallest positive number where either $x(t_1)=0$ or $y(t_1)=0$. 
However, the above argument implies $$x(t_1) \geq x_0 \exp\left\{ \int_{t_0}^{t_1} r_1(s)~ds \right\}>0,~~
y(t_1) \geq y_0 \exp\left\{ \int_{t_0}^{t_1} r_2(s)~ds \right\}>0,$$ which is  a contradiction, hence
all solutions of \eqref{1a},\eqref{2star} are positive. This also excludes the possibility
that either (\ref{inf1}) or a similar equality for $y(t)$ holds and concludes the proof of Part 1) in the statement 
of the theorem.

2) Assuming (a2$^{\ast}$), let us prove that (\ref{inf2}) cannot be satisfied. By the assumption in (a6),
both initial functions are bounded.
Fix some $\varepsilon>0$ and denote 
$\nu_1=\max\{ K+\varepsilon,
\sup_{s \leq 0} \varphi(s)+\varepsilon\}$, 
$\nu_2=\max\{ f_2(K)+\varepsilon,\sup_{s \leq 0} \psi(s)+\varepsilon\}$. 
Let us verify that there exist positive bounds $M_1$, $M_2$ for the solutions 
$x$ and $y$, respectively, such that $f_2(M_1) < M_2 < f_1^{-1}(M_1)$
and $f_1(M_2) < M_1 < f_2^{-1}(M_2)$, which means that the point $(M_1,M_2)$ is
between the curves $f_2(x)$ (the lower curve) and $f_1^{-1}(x)$ (the upper curve), $x>K$.

If $f_2^{-1}(\nu_2) \leq \nu_1$, denote $M_2=f_2(\nu_1)+\varepsilon_1$, where
$\varepsilon_1< f_1^{-1}(\nu_1)-f_2(\nu_1)$ is a positive number, which exists 
since $f_1^{-1}(x)-f_2(x)>0$ for $x>K+\varepsilon$, and $M_1=\nu_1$.

If $f_2^{-1}(\nu_2) >\nu_1$ but $f_1^{-1}(\nu_1) \leq \nu_2$, denote $M_1=f_2^{-1}(\nu_2)+\varepsilon_1$, where
$\varepsilon_1< f_1(\nu_2)-f_2^{-1}(\nu_2)$, and $M_2=\nu_2$.

If both $f_2^{-1}(\nu_2) >\nu_1$ and $f_1^{-1}(\nu_1) >\nu_2$, we can take $M_1=\nu_1$ and $M_2=\nu_2$,
then $f_2(M_1) < M_2 < f_1^{-1}(M_1)$ and $f_1(M_2) < M_1 < f_2^{-1}(M_2)$.

We have $x_0<M_1$, $y_0<M_2$, these inequalities are also valid on $[0, t_0)$ for some
$t_0>0$. Let us prove that $x(t) < M_1$, $y(t)<M_2$ for any $t\geq 0$. Let us assume the
contrary, and let $t_1$ be the smallest point where either $x(t_1)=M_1$ or $y(t_1)=M_2$.
Suppose $x(t_1)=M_1$, the case $y(t_1)=M_2$ is considered similarly.    
Denote $t^{\ast}=\sup\{ t \in [0,t_1] | x(t) \leq f_1(M_2) \}$, so $x(t^{\ast})=f_1(M_2)$
and $x(t)>f_1(M_2)$ for $t \in (t^{\ast}, t_1]$
However, for $t \in [t^{\ast},t_1]$ we have $y(t) \leq M_2$, so $f_1(y(t))\leq f_1(M_2)$ and
$f_1(M_2)<x(t)<M_1$, thus due to monotonicity of $f_1$
$$
\frac{dx}{dt}  = r_1(t) \left[ \int_{h_1(t)}^t f_1(y(s))~d_s R_1 (t,s) - x(t) \right]
\leq r_1(t) \left[ f_1(M_2)-f_1(M_2) \right] =0,$$
non-positivity of the derivative of $x$ on $[t^{\ast},t_1]$ implies $M_1=x(t_1)\leq x(t^{\ast})=f_1(M_2)$, which is a 
contradiction. Thus (\ref{inf2}) is impossible and there exists a positive global solution.

3) Next, assume that (a2) holds, which is a particular case of (a2$^{\ast}$), 
and prove permanence of equation (\ref{1a}) with positive initial conditions.

By (a6) we have $x_0=x(0)>0$, $y_0=y(0)>0$, and according to (a3) there is $t_1$
such that $h_i(t) \geq 0$ for $t\geq t_1$, $i=1,2$. 
From positivity of solutions justified in Part 2, there are $\mu_1$ and $\mu_2$ such that
$x(t) \geq\mu_1>0$ and $y(t) \geq \mu_2>0$ for any $t \in [0,t_1]$. 
In particular, we can choose $m_1>0$ and $m_2>0$ satisfying 
\begin{equation}
\label{bound0}
m_1 < \min\left\{ \mu_1, f_2^{-1}(\mu_2),K  \right\}, ~~
m_2 < \min\left\{ \mu_2, f_1^{-1}(\mu_1),K  \right\},
\end{equation}
and also such that the point $(m_1,m_2)$ is between the curves $y=f_1^{-1}(x)$ and 
$y=f_2(x)$, where $0<m_1<K$, so 
\begin{equation}
\label{ball1}
f_1(m_2)>m_1, ~~ f_2(m_1)>m_2,
\end{equation}
which is possible since $f_2(x)>f_1^{-1}(x)$ for 
$x\in [0,K]$ and thus in $(0, \min\left\{ \mu_1, f_2^{-1}(\mu_2),K  \right\})$.
Thus any point $(m_1,m_2)$ between the curves $y=f_1^{-1}(x)$ and $y=f_2(x)$,  
see Fig.~\ref{figure1}, satisfies \eqref{ball1}.

Further, let us verify that $x(t) \geq m_1$, $y(t) \geq m_2$ for any $t \geq 0$. 
As defined, $x(t)\geq m_1$, $y(t) \geq m_2$ for $t \in [0,t_1]$, and also $h_i(t) \geq 0$ 
for $t\geq t_1$, $i=1,2$. Thus $x(t)$ is greater than the solution of the ordinary differential equation 
$$
x'(t)=r_1(t)[f_1(m_2)-x(t)]
$$
as long as $y(t) \geq m_2$, so $x(t)$ is increasing if $m_1<x(t)<f_1(m_2)$, thus $x(t)>m_1$ unless $y(t)$ becomes smaller than $m_2$ (in fact, even smaller
than $f^{-1}(m_1)<m_2$). However, 
$$
y'(t)>r_2(t)[f_2(m_1)-y(t)]
$$
as long as $x(t) \geq m_1$, thus $x(t) \geq m_1$, $y(t) \geq m_2$ for any $t \geq 0$.

The upper bound of the solution was constructed in Part 2), thus the solution is permanent, which concludes the proof.
 
\begin{corol}
\label{th1_cor1}
The results of Theorem~\ref{theorem1} hold for system \eqref{1avar},\eqref{2star} if instead of (a4)
we assume (b4). 
\end{corol}

\begin{corol}
\label{th1_cor2}
The results of Theorem~\ref{theorem1} hold for system \eqref{1aintegro},\eqref{2star} if assumption (a4) is replaced by 
(c4). 
\end{corol}

\begin{remark}
Let us note that (a2$^{\ast}$) guarantees global boundedness but not persistence of solutions,
see Example~\ref{example2} where the solution tends to (0,0) as $t \to \infty$.
\end{remark}

The following examples illustrate the fact that when (a2) is not satisfied, 
the solution can fail to be either bounded or persistent, even for a non-delay system.

\begin{example}
\label{example1}
Let $f_i(x)=x^2+x$, $r_i(t)=1$, $h_i(t)=t$ in \eqref{1avar}. The system
\begin{equation}
\label{ex1eq1}
x'(t)=y^2(t)+y(t)-x(t), ~~
y'(t)=x^2(t)+x(t)-y(t)
\end{equation}
has an unbounded solution $\displaystyle (x(t),y(t))= \left( \frac{1}{3-t}, \frac{1}{3-t}  \right)$ on $[0,3)$. 
The functions $f_i(x)$ satisfy $f_2(x)>x>f_1^{-1}(x)$, so (a2) does not hold, there is no positive equilibrium.
\end{example}

\begin{example}
\label{example2}
For $f_i(x)=\frac{1}{2}x$, $r_i(t)=1$, $h_i(t)=t$ in \eqref{1avar}, the system
\begin{equation}
\label{ex2eq1}
x'(t)=\frac{1}{2}y(t)-x(t), ~~
y'(t)=\frac{1}{2}x(t)-y(t)
\end{equation}
has a solution $\displaystyle (x(t),y(t))= \left( e^{-t/2},  e^{-t/2} \right)$ on $[0,\infty)$ 
which tends to zero as $t \to \infty$ and thus is not persistent. 
The functions $f_1(x)=f_2(x)= \frac{1}{2}x$ satisfy $f_2(x)=\frac{1}{2}x<2x=f_1^{-1} (x)$ for any 
$x>0$, thus (a2$^{\ast}$) is satisfied while (a2) is not.
\end{example}

\section{Stability of the Positive Equilibrium}
\label{stability}

Next, let us proceed to stability.
The following result considers the case when $f_2(x)>f_1^{-1}(x)$ on $(0,K)$, and $f_2(x)<f_1^{-1}(x)$ on $(K,\infty)$.
This can be interpreted as cooperation for small $x$ and competition for large $x$. Theorem~\ref{theorem2} states
that in this case the equilibrium $(K,f_2(K))$ attracts all positive solutions.

\begin{figure}[ht]
\centering
\includegraphics[scale=0.95]{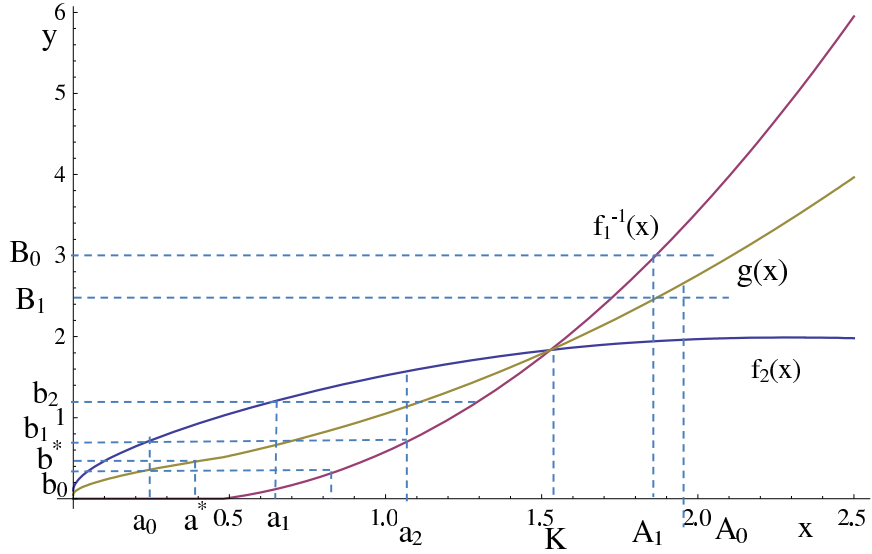} 
\caption{Convergence and solution bounds}
\label{figure2}   
\end{figure}

\begin{guess}
\label{theorem2}
Suppose (a1)-(a6) hold. Then any solution of \eqref{1a},\eqref{2star} converges to the unique positive equilibrium 
$(x(t),y(t)) \to (K,f_2(K))$ as $t \to \infty$.
\end{guess}
{\bf Proof.} The proof is illustrated by Fig.~\ref{figure2}.

According to Theorem~\ref{theorem1}, there are $a,A,b,B$ such that $0<a \leq x(t) \leq A$ and $0<b \leq y(t) \leq B$
for any $t \geq 0$. We can always assume $a<K<A$ and $b<K<B$ without loss of generality.

Consider in addition to $f_1$, $f_2$ a monotone increasing function $g:\RR^+ \to \RR^+$ 
which satisfies $f_1^{-1}(x)<g(x)<f_2(x)$ on $(0,K)$ and $f_2(x)<g(x)<f_1^{-1}(x)$ on $(K,\infty)$; in particular, we can 
take $g(x)= \alpha f_1^{-1}(x)+ (1-\alpha)f_2(x)$, $\alpha \in (0,1)$, where we assume $f_1^{-1}(x)=0$
if there is no non-negative $t$ such that $f_1(t)=x$.
As we assumed $f_1^{-1}(0)=0$, we can always find $\alpha \in (0,1)$ such that 
$g(0) \leq b$.

Further, let us choose $a_0=\min\{ a, g^{-1}(b) \}$, $b_0=\min\{  b, g(a) \}$. The function
$g(x)$ is monotone increasing, and we have either $g(a)\geq b$ or $g(a)<b$.
In the former case $b_0=b$ and $a>g^{-1}(b)$, so $a_0=g^{-1}(b)$ and $b_0=g(a_0)$. In the latter case 
$b_0=g(a)$ and $a<g^{-1}(b)$, so $a_0=a$ and $b_0=g(a_0)$.

We also have $a_0\leq a$, $b_0 \leq b$, so 
\begin{equation}
\label{0iteration_first}
0<a_0 \leq x(t) \leq A_0, ~~~0<b_0 \leq y(t) \leq B_0, ~~t \geq 0.
\end{equation}
By (a3), there is $t_0>0$ such that $h_i(t)\geq 0$ for $t \geq t_0$, $i=1,2$.
Next, we have $f_1(b_0)>b_0$, $f_2(a_0)>a_0$ since $g(a_0)=b_0$, and 
the curve $g(x)$ is between $f_1^{-1}(x)$ and $f_2(x)$, see Fig.~\ref{figure2}.

Define $a_1=\min\{ g^{-1}(f_2(a_0)), f_1(b_0)  \}$, $b_1=\min\{ f_2(a_0), g(f_1(b_0)) \}$,
where 
\begin{equation}
\label{prep_1_iteration}
\begin{array}{lll}
b_1=g(a_1), & \displaystyle f_1(y) \geq f_1(b_0)-a_1 \geq 0 & \mbox{~~for~~} y \in [b_0,b_1], \\
 & \displaystyle f_2(x) \geq f_1^{-1}(a_0)-b_1 \geq 0 & \mbox{~~for~~} x\in [a_0,a_1],
\end{array}
\end{equation}
and the inequalities are strict for any $x<a_1$, $y<b_1$. 
Thus \eqref{1a} and \eqref{prep_1_iteration}  imply $x'(t)>0$, $y'(t)>0$ for $(x,y) \in [a_0, a_1]\times [b_0,b_1]$,
and the derivative is positive for any $x<a_1$, $y<b_1$. Let us prove that there exists $t_1 \geq t_0$ such that
\begin{equation}
\label{1iteration}
0<a_1 \leq x(t), ~~~0<b_1 \leq y(t), ~~t \geq  t_1.
\end{equation}

Let us choose $a^{\ast} \in (a_0,a_1)$, $b^{\ast} \in (b_0,b_1)$ and first prove that there exists $t^{\ast}$ such that
$x(t^{\ast}) \geq a^{\ast}$, $y(t^{\ast}) \geq b^{\ast}$. If $x(t_0)$ and $y(t_0)$ satisfy these inequalities, they are also satisfied 
for any $t \geq t_0$ due to \eqref{prep_1_iteration}, and there is nothing to prove. If either $x(t_0)< a^{\ast}$ or $y(t_0)< 
b^{\ast}$, or both, 
then the derivative exceeds a positive value
$$x'(t)> r_1(t)[ f_1(b_0)-a^{\ast}], ~~y'(t) >r_2(t) [ f_1^{-1}(a_0)-b^{\ast}],$$
as long as $x(t)<a^{\ast}$, $y(t) <b^{\ast}$,
where the expressions in the brackets are positive constants. Due to (a5), 
there is a point $t^{\ast}$ such that $x(t^{\ast}) \geq a^{\ast}$, $y(t^{\ast}) \geq b^{\ast}$.
Moreover, as \eqref{0iteration_first} holds, these inequalities are satisfied for $t \geq t^{\ast}$ as well. 
Let us choose $\bar{t}$ such that 
$h_i(t)\geq t^{\ast}$ for $t \geq \bar{t}$, $i=1,2$. Then
$$x'(t)> r_1(t)[ f_1(b^{\ast})-a_1], ~~y'(t) >r_2(t) [ f_1^{-1}(a^{\ast})-b_1],~~t \geq \bar{t},$$
as long as $x(t)<a_1$, $y(t) <b_1$, and the expressions in the brackets are positive constants.
Again, referring to (a5), we obtain that there exists $t_1^{\ast} \geq \bar{t}$ such that \eqref{1iteration} holds.
Applying the same procedure to the upper bound, we find $t_1 \geq t_1^{\ast}$
such that
\begin{equation}
\label{1iter}
0<a_1 \leq x(t) \leq A_1, ~~~ 0<b_1 \leq y(t) \leq B_1, ~~t \geq t_1.
\end{equation}
Continuing the process by induction, we obtain increasing sequences $\{ a_n\}$, $\{ b_n \}$, $\{ t_n \}$
and decreasing sequences $\{ A_n\}$, $\{ B_n \}$, where $g(a_n)=b_n$, $g(A_n)=B_n$ and
\begin{equation}
\label{0iteration}
0<a_n \leq x(t) \leq A_n, ~~~0<b_n \leq y(t) \leq B_n, ~~t \geq t_n.
\end{equation}
Thus all the sequences have limits: $\lim_{n \to \infty} a_n=a$, $\lim_{n \to \infty} b_n=b$, and $g(a)=b$; moreover, all $a_n < K$, $b_n<K$, so $a \leq K$,
$b \leq K$. If $a<K$ then $f_1(b)>a$ and $f_2(a)>b$, and from continuity there exists $\varepsilon >0$ such that 
$f_2(x)>b$ for $x \in (a-\varepsilon, a)$ and $f_1(y)>a$ for $y \in (g^{-1}(a-\varepsilon),b)$.
As $a$ is a limit, there exists $a_k\in (a-\varepsilon, a)$, then $a_{k+1} = \min\{ g^{-1}(f_2(a_k)), f_1(b_k)  \}>a$,
which leads to a contradiction with $a>a_k$ for any $k$. Hence $a=K$; similarly, we can prove that $A=K$ and thus  
any solution of \eqref{1a},\eqref{2star} converges to the unique positive equilibrium: 
$(x(t),y(t)) \to (K,f_2(K))$ as $t \to \infty$.

\begin{corol}
\label{th2_cor1}
The results of Theorem~\ref{theorem2} hold for system \eqref{1aintegro},\eqref{2star} if assumption (a4) is replaced by 
(b4). 
\end{corol}

\begin{corol}   
\label{th2_cor2}
The results of Theorem~\ref{theorem2} hold for system \eqref{1aintegro},\eqref{2star} 
if assumption (a4) is replaced by (c4).
\end{corol}

\begin{example}
\label{example2a}
Let us note that condition (a5) is not required for permanence of solutions but is crucial for convergence to 
the unique positive equilibrium. 
For example, if $f_i(x)=1+\frac{1}{2}x$, $r_1(t)=2+\sin t$, 
$r_2(t)=2+\cos t$, $K_i(t,s)=1$, $h_i(t)=t-1$ in \eqref{1aintegro}, then
 all solutions with positive initial values and non-negative initial functions of the system
\begin{equation}
\label{ex2aeq1}
\begin{array}{lll}
x'(t) & = & \displaystyle (2+\sin t) \left[ \int_{t-1}^t \left( 1+\frac{1}{2} y(s) \right)~ds - x(t) \right], \vspace{2mm} 
\\
y'(t) & = & \displaystyle (2+\cos t)  \left[ \int_{t-1}^t  \left[ 1+\frac{1}{2} x(s) \right)~ds -y(t) \right]
\end{array}
\end{equation}
converge to the unique positive equilibrium point (2,2), since all the conditions of Theorem~\ref{theorem2}
are satisfied.

However, system \eqref{1avar} with $f_i(x)=1+\frac{1}{2}x$, $\displaystyle r_i(t)=\frac{2}{e^{2t}+0.5}$,
$h_i(t)=t$, which is 
\begin{equation}
\label{ex2aeq2}
\begin{array}{l}
x'(t)= \displaystyle \frac{2}{e^{2t}+0.5} \left[ \left( 1+\frac{1}{2} y(t) \right)-x(t)\right], \vspace{2mm} \\
y'(t)= \displaystyle \frac{2}{e^{2t}+0.5} \left[ \left(  1+\frac{1}{2} x(t)\right)-y(t)\right],
\end{array}
\end{equation}
has a solution $\displaystyle \left( 4+e^{-2t}, 4+e^{-2t} \right)$ which tends to $(4,4)$ as $t \to \infty$,
not to the unique positive equilibrium point (2,2). For system \eqref{ex2aeq2} with $(x(0),y(0))=(5,5)$,
all the conditions of Theorem~\ref{theorem2} but (a5) are satisfied, since $\displaystyle \int_0^{\infty}
r(t)~dt = \int_0^{\infty} \frac{2~dt}{e^{2t}+0.5}<\infty$.
\end{example}

In contrast to Theorem~\ref{theorem2}, if $f_2(x)<f_1^{-1}(x)$ for any $x$, all positive solutions converge to zero,
which can be interpreted as a continuing negative mutual influence leading to extinction.
In the case $f_2(x)>f_1^{-1}(x)$ for any $x$, all positive solutions are unbounded and tend to infinity.
The effect is due to mutual positive feedback.

\begin{figure}[ht]
\centering
\includegraphics[scale=0.9]{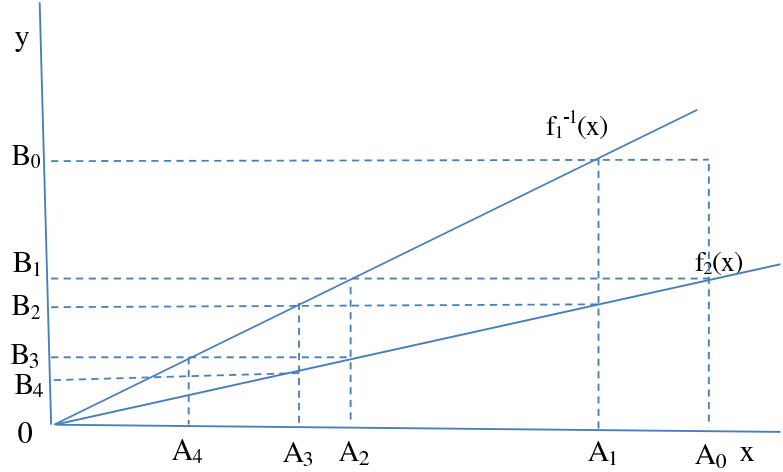} 
\caption{Illustration of the solution bounds tending to zero}
\label{figure3}   
\end{figure}

\begin{guess}
\label{theorem3}
Suppose (a1) and (a3)-(a6) hold. 

1) If $f_2(x)<f_1^{-1}(x)$ for any $x>0$ then every solution of \eqref{1a},\eqref{2star} converges to zero as $t \to 
\infty$.

2) If $f_2(x)>f_1^{-1}(x)$ for any $x>0$ then every global solution of \eqref{1a},\eqref{2star} 
tends to $+\infty$ as $t \to \infty$.
\end{guess}
{\bf Proof.} 1) The proof is similar to the proof of Theorem~\ref{theorem2}. 
First, we notice that there exist $A_0,B_0>0$ such that
$$
0< x(t) \leq A_0, ~~~ 0< y(t) \leq B_0, ~~t \geq t_0.
$$
Next, we define $A_1=f_1(B_0)$, $B_1=f_2(A_0)$ (see Fig.~\ref{figure3}) and prove that for some $t_1>t_0$ we have
$$
0< x(t) \leq A_1, ~~~ 0< y(t) \leq B_1, ~~t \geq t_1.
$$

By induction we verify
$$
0< x(t) \leq A_n, ~~~ 0< y(t) \leq B_n, ~~t \geq t_n
$$
where $A_n=f_1(B_{n-1})$, $B_n=f_2(A_{n-1})$ and both sequences $\{A_n\}$ and $\{ B_n \}$ are positive, decreasing 
and hence have a limit. Let $\displaystyle d = \lim_{n \to \infty} A_n$, then by construction and continuity of $f_i$ we have
$\displaystyle \lim_{n \to \infty} B_n =f_2(d)$ and $f_1(f_2(d))=d$, so $d=f_2(d)=0$. Thus any solution of \eqref{1a},\eqref{2star} converges to zero as $t \to \infty$.

The proof of 2) is similar.

\begin{corol}
\label{th3_cor1} 
The results of Theorem~\ref{theorem3} hold for system \eqref{1avar},\eqref{2star} if assumption (a4) is replaced by
(b4).
\end{corol}

\begin{corol}   
\label{th3_cor2}
The results of Theorem~\ref{theorem3} hold for system \eqref{1aintegro},\eqref{2star} if assumption (a4) is replaced by
(c4).
\end{corol}

\begin{example}
\label{example3}
For any $h>0$, consider system
\eqref{1aintegro} with $f_i(x)=x^2+x$, $\displaystyle r_i(t)=1$,
$h_i(t)=t-h$, $\displaystyle K_i(t,s)=\frac{1}{h}$, with
\begin{equation}
\label{ex3eq1}
x'(t)= \frac{1}{h} \int_{t-h}^t \left[ y^2(s)+y(s)\right]~ds-x(t), ~~
y'(t)= \frac{1}{h} \int_{t-h}^t \left[ x^2(s)+x(s)\right]~ds-y(t).
\end{equation}
Every solution with non-negative initial conditions and positive initial values tends to infinity at the right end of the 
maximal interval where the solution exists, which illustrates Part 2 of Theorem~\ref{theorem3}.
\end{example}

\begin{example}
\label{example4}
For any $h>0$, consider system
\eqref{1aintegro} with $\displaystyle f_1(x)=e^x-1$, $\displaystyle f_2(x)=\frac{1}{2} \ln(x+1)$, $\displaystyle r_i(t)=1$,
$h_i(t)=t-h$, $\displaystyle K_i(t,s)=\frac{2}{h^2}(s+h-t)$, which is
\begin{equation}
\label{ex4eq1}
\begin{array}{l}
\displaystyle x'(t)=\frac{2}{h^2}  \int_{t-h}^t (s+h-t) \left( e^{y(s)}-1 \right)~ds -x(t), \vspace{2mm} \\ 
y'(t)=\frac{2}{h^2}  \int_{t-h}^t (s+h-t) \frac{1}{2} \ln(x(s)+1)~ds-y(t).
\end{array}
\end{equation}
Then $f_1^{-1}(x)=\ln(x+1)>f_2(x)$ for positive $x$, by Part 1of Theorem~\ref{theorem3} every solution with 
non-negative initial conditions and positive initial values tends to zero as $t \to \infty$.
\end{example}

In both Theorems~\ref{theorem2} and ~\ref{theorem3}, all positive solutions had the same asymptotics.
Theorem~\ref{theorem4} considers the case when the limit behaviour depends on the initial conditions.
In particular, two cases are considered. In the first case,
for small initial conditions, a solution tends to zero, while
for large initial conditions, a solution tends to the unique positive
equilibrium.
In the second case,  
for small initial conditions, a solution tends to the unique positive 
equilibrium, for large initial conditions, a solution tends to infinity.

\begin{guess}
\label{theorem4}
Suppose (a1) and (a3)-(a6) hold. 

1) If $f_2(x)<f_1^{-1}(x)$ 
for any $x>0$, $x\neq K$ and   $f_2(K)=f_1^{-1}(K)$
then any solution of \eqref{1a},\eqref{2star} with the initial function satisfying $\varphi(t)\geq K$,
$\psi(t) \geq f_2(K)$ tends to $(K,f_2(K))$ as $t \to \infty$, while any solution of \eqref{1a},\eqref{2star} 
with the initial function satisfying $\varphi(t)< K$,
$\psi(t) < f_2(K)$ converges to zero as $t \to \infty$.

2) If $f_2(x)>f_1^{-1}(x)$ for any $x>0$, $x\neq K$  and $f_2(K)=f_1^{-1}(K)$
then any solution of \eqref{1a},\eqref{2star} with the initial function satisfying $\varphi(t)\geq K$,
$\psi(t) \geq f_2(K)$ tends to $+\infty$ as $t \to \infty$ 
(or $t \to c^-$, where $c$ is a finite right end of the maximal interval of the existence of the solution), while any 
solution of \eqref{1a},\eqref{2star} with the initial function satisfying $\varphi(t)< K$,
$\psi(t) < f_2(K)$ converges to $(K,f_2(K))$ as $t \to \infty$.
\end{guess} 
{\bf Proof.} 1) The proof of the case $\varphi(t)\geq K$, $\psi(t) \geq f_2(K)$ completely coincides with the proof of the 
upper bound in Theorem~\ref{theorem2}, since the lower bound of the solution is $K$ and, as in (a2), $f_2(x)<f_1^{-1}(x)$ for $x>K$.
If $\varphi(t)\in (0,K)$, $\psi(t) \in (0,f_2(K))$, then we repeat the previous proof for $f_2(x)<f_1^{-1}(x)$,
where the zero takes the place of $K$.

Similarly, in 2) the part $\varphi(t)\in (0,K)$, $\psi(t) \in (0,f_2(K))$ coincides with the proof
that the lower bound tends to $K$ in Theorem~\ref{theorem2}, as $f_2(x)>f_1^{-1}(x)$, $x \in (0,K)$.
For the proof of the second part we construct a sequence of upper bounds which tends to $+\infty$,
as in the proof of Theorem~\ref{theorem3}.

\begin{corol}
\label{th4_cor1}
The results of Theorem~\ref{theorem4} hold for system \eqref{1avar},\eqref{2star} if assumption (a4) is replaced by
(b4).
\end{corol}

\begin{corol}
\label{th4_cor2}
The results of Theorem~\ref{theorem4} hold for system \eqref{1aintegro},\eqref{2star} if assumption (a4) is replaced by
(c4).
\end{corol}

The following result can be interpreted as nonoscillation about the unique positive equilibrium.

\begin{guess}
\label{theorem5}
Suppose (a1)-(a4) and (a6) hold.
Any solution of \eqref{1a},\eqref{2star} with the initial function satisfying $\varphi(t)\geq K$,
$\psi(t) \geq f_2(K)$ satisfies $x(t)\geq K$, $y(t) \geq f_2(K)$ for any $t \geq 0$, while 
any solution of \eqref{1a},\eqref{2star} with $\varphi(t)\leq K$, $\psi(t) \leq f_2(K)$ satisfies $x(t)\leq K$, $y(t) \leq f_2(K)$ for 
any $t \geq 0$.
\end{guess}
{\bf Proof.} Consider the case  $\varphi(t)\geq K$, $\psi(t) \geq f_2(K)$. 
As long as $x(t)\geq K$, $y(t) \geq f_2(K)$, we have 
$$
\frac{dx}{dt} \geq  r_1(t) \left[ \int_{h_1(t)}^t f_1(f_2(K))~d_s R_1 (t,s) - x(t) \right]=
r_1(t) [K-x(t) ], $$ $$  
\frac{dy}{dt} \geq r_2(t) \left[ \int_{h_2(t)}^t f_2(K)~d_s R_2 (t,s) - y(t) \right] =r_2(t)  [f_2(K) - y(t)].
$$
The first inequality implies $\displaystyle x(t) \geq \varphi(0)-K+K \exp\left\{ \int_0^t r_1(s)~ds \right\}$
which exceeds $K$ for $\varphi(0)>K$ and is identically equal to $K$ if $\varphi(0)=K$.
The second inequality gives $\displaystyle y(t) \geq \psi(0)-f_2(K)+f_2(K) \exp\left\{ \int_0^t r_2(s)~ds \right\}$
which is also not less than $f_2(K)$. 

Again, using monotonicity of $f_i$, the case $0 \leq \varphi(t)\leq K$, $0 \leq \psi(t) \leq f_2(K)$,
$\varphi(0)>0$, $\psi(0)>0$ is treated in a similar way.

A more general model
\begin{equation}
\label{1b}
\begin{array}{ll}
\displaystyle \frac{dx}{dt} & \displaystyle = r_1(t)G_1(x) \left[ \int_{h_1(t)}^t f_1(y(s))~d_s R_1 (t,s) - x(t) \right] 
\vspace{2mm} \\
\displaystyle \frac{dy}{dt} & \displaystyle = r_2(t)G_2(y) \left[ \int_{h_2(t)}^t f_2(x(s))~d_s R_2 (t,s) - y(t) \right]
\end{array}
\end{equation}
includes a system of logistic equations with the delay in the production term described in \cite{AWW}.
We assume that the functions $G_i$ satisfy 
\begin{description}
\item{{\bf (a7)}} 
$G_i:\RR^+ \to \RR^+$, $i=1,2$ are continuous functions satisfying 
$G_i(x)>0$ for $x>0$.
\end{description}

The proofs of the following results coincide with the proofs of 
Theorems~\ref{theorem1},\ref{theorem2},\ref{theorem3},\ref{theorem4}, respectively.

\begin{guess}
\label{theorem1a}
Suppose (a1)-(a4),(a6)-(a7) hold. Then any solution of \eqref{1b},\eqref{2star} is permanent.
\end{guess}

\begin{guess}   
\label{theorem2a}
Suppose (a1)-(a7) hold. Then any solution of \eqref{1b},\eqref{2star} converges to the unique positive equilibrium
$(x(t),y(t)) \to (K,f_2(K))$ as $t \to \infty$.
\end{guess}   

\begin{guess}
\label{theorem3a}
Suppose (a1) and (a3)-(a7) hold.

1) If $f_2(x)<f_1^{-1}(x)$ for any $x>0$ then any solution of \eqref{1b},\eqref{2star} converges to zero as $t \to \infty$.

2) If $f_2(x)>f_1^{-1}(x)$ for any $x>0$ then any global solution of \eqref{1b},\eqref{2star}
tends to infinity as $t \to \infty$.
\end{guess}

\begin{guess}
\label{theorem4a}
Suppose (a1) and (a3)-(a7) hold.

1) If $f_2(x)<f_1^{-1}(x)$
for any $x>0$, $x\neq K$ and   $f_2(K)=f_1^{-1}(K)$  
then any solution of \eqref{1b},\eqref{2star} with the initial function satisfying $\varphi(t)\geq K$,
$\psi(t) \geq f_2(K)$ tends to $(K,f_2(K))$ as $t \to \infty$, while any solution of \eqref{1b},\eqref{2star}
with the initial function satisfying $\varphi(t)< K$,
$\psi(t) < f_2(K)$ converges to zero as $t \to \infty$.

2) If $f_2(x)>f_1^{-1}(x)$ for any $x>0$, $x\neq K$  and $f_2(K)=f_1^{-1}(K)$
then any solution of \eqref{1b},\eqref{2star} with the initial function satisfying $\varphi(t)\geq K$,
$\psi(t) \geq f_2(K)$ tends to $+\infty$ as $t \to \infty$ 
(or $t \to c^-$, where $c$ is a finite right end of the maximal interval of the existence of the solution),
while any solution of \eqref{1b},\eqref{2star} with $\psi(t) < f_2(K)$ converges to $(K,f_2(K))$ as $t \to \infty$.
\end{guess}

\begin{guess}
\label{theorem5a} 
Suppose (a1)-(a4) and (a6)-(a7) hold.
Any solution of \eqref{1a},\eqref{2star} with the initial function satisfying $\varphi(t)\geq K$,
$\psi(t) \geq f_2(K)$ satisfies $x(t)\geq K$, $y(t) \geq f_2(K)$ for any $t \geq 0$, while
any solution of \eqref{1a},\eqref{2star} with $\varphi(t)\leq K$, $\psi(t) \leq f_2(K)$ 
satisfies $x(t)\leq K$, $y(t) \leq f_2(K)$ for any $t \geq 0$.
\end{guess}

\begin{example}
\label{example5}
All solutions of system \eqref{1a} with $G_i(x)=x$, $r_i(x)=1$,
$R_i(t,s)=\chi_{(h_i(t),\infty)}(s)$, where $\chi_I(\cdot)$ is the characteristic function of set $I$,  
$$
x'(t)=x(t)\left[ \sqrt{y(h_1(t))}+2-x(t) \right],~~
y'(t)=y(t)\left[ x(h_2(t))-y(t) \right],$$
with non-negative initial functions and non-trivial initial value
converge to the positive equilibrium (4,4). This is also true for solutions
with non-negative initial functions and non-trivial initial value 
of  system \eqref{1a} with $f_1(x)=\sqrt{x}+2$, $f_2(x)=x$, $\displaystyle r_i(t)=1$,
$h_i(t)=t-h$, $R_i(t,s)=\frac{1}{h}(s-t+h)$ if $s \in [t-h,t]$ and zero elsewhere,
which is
$$
x'(t)= \frac{x(t)}{h} \int_{t-h}^t (\sqrt{y(s)}+2)~ds-x^2(t),~~
y'(t)= \frac{y(t)}{h} \int_{t-h}^t x(s)~ds-y^2(t).$$
\end{example}

\section{Applications and Discussion}
\label{summary}

As an example, consider the following models of type (\ref{intro6}) 
\begin{equation}
\label{summary1}
\begin{array}{ll}
\displaystyle \frac{dx}{dt} & \displaystyle = c_1 \int_{h_1(t)}^t K_1(t,s) \tanh(y(s))~ds - \mu_1 x(t)   \vspace{2mm} \\
\displaystyle \frac{dy}{dt} & \displaystyle = c_2 \int_{h_2(t)}^t K_2(t,s) \tanh(x(s))~ds - \mu_2 y(t)
\end{array}
\end{equation}
and
\begin{equation}
\label{summary2}
\frac{dx}{dt} = c_1 \tanh(y(h_1(t))-\mu_1 x(t), ~~ \frac{dy}{dt}  = c_2 \tanh(x(h_2(t))-\mu_2 y(t),
\end{equation}
where $c_1,c_2,\mu_1,\mu_2$ are positive constants, $\displaystyle \int_{h_1(t)}^{t} K_1(t,s)~ds=\int_{h_2(t)}^t K_2(t,s)~ds \equiv 
1$, $K_i(t,s) \geq 0$ for any $t\geq 0$, $s>0$, $h_i$ are Lebesgue measurable functions satisfying
$h_i(t)\geq 0$ and $\displaystyle \lim_{t \to \infty} h_i(t)=\infty$, $i=1,2$.
A more general version of \eqref{summary1} but with constant delays was studied in \cite{Dong}.

Since for $u>0$ the functions $f_i(u)=c_i \tanh(u)/\mu_i$ satisfy $f_i^{\prime}(u)>0$, $f_i^{\prime\prime}(u)<0$,
there is a positive equilibrium for $f_2^{\prime}(0)>1/f_1^{\prime}(0)$, or $c_2/\mu_2>\mu_1/c_1$, otherwise
$f_2(x)<f_1^{-1}(x)$ for $x>0$.

Theorems \ref{theorem2} and \ref{theorem3} imply the following result.

If $c_1c_2 >\mu_1 \mu_2$, then all solutions of \eqref{summary1},\eqref{2star} and \eqref{summary2},\eqref{2star}  with 
non-negative initial functions and positive initial values converge to the unique positive equilibrium
$(x^{\ast},y^{\ast} )$, where $x^{\ast}$ is a solution of the 
equation
$$ \frac{c_1}{\mu_1} \tanh \left( \frac{c_2}{\mu_2} \tanh(x^{\ast}) \right)=x^{\ast},$$
$y^{\ast}=c_2\tanh(x^{\ast})/\mu_2$.
If $c_1c_2 \leq \mu_1 \mu_2$, then all solutions of \eqref{summary1},\eqref{2star} and \eqref{summary2},\eqref{2star}
converge to zero.

Next, for the Lotka-Volterra-type cooperative system
\begin{equation}
\label{summary3}
\begin{array}{ll}
\displaystyle \frac{dx}{dt} & \displaystyle = r_1(t) x(t) \left[ A_1 - a_1 x(t) + b_{1} \int_{h_1(t)}^t y(s)~d_s R_1(t,s) 
\right]  \vspace{2mm} \\
\displaystyle \frac{dy}{dt} & \displaystyle  = r_2(t) y(t) \left[ A_2 - a_2 y(t) + b_{2} \int_{h_2(t)}^t x(s)~d_s R_2(t,s)  
\right]
\end{array}
\end{equation} 
Theorems~\ref{theorem2} and \ref{theorem3} imply the following result.

\begin{guess}
\label{theorem_appl1}
Suppose (a3)-(a6) hold, $A_i\geq 0$, $a_i>0$ and  $b_i>0$, $i=1,2$.  

If $A_1+A_2>0$ and $a_1a_2<b_1b_2$ then there exists a unique positive equilibrium
$((b_1A_2-a_2A_1)/(b_1b_2-a_1a_2), (b_2A_1-a_1A_2)/(b_1b_2-a_1a_2))$ of \eqref{summary3}, and
all solutions of \eqref{summary3},\eqref{2star} converge to this equilibrium.
If $A_1=A_2=0$ and $a_1a_2<b_1b_2$ all solutions of \eqref{summary3},\eqref{2star} 
converge to (0,0). If $A_1+A_2>0$ and $a_1a_2 \geq b_1b_2$ then both components of the solution of 
\eqref{summary3},\eqref{2star} tend to $+\infty$ as $t \to \infty$.
\end{guess}

Next, let us consider the generalization of the cooperative model \cite[p.192]{Gopalsamy}
to the case of distributed delays and time-variable growth rates
\begin{equation}
\label{summary4}
\begin{array}{ll}
\displaystyle \frac{dx}{dt} & 
\displaystyle = r_1(t) x(t) \left[  \int_{h_1(t)}^t \frac{K_1+\alpha_1 y(s)}{1+y(s)}~d_s R_1(t,s) - x(t) \right]  
\vspace{2mm} \\
\displaystyle \frac{dy}{dt} & 
\displaystyle = r_2(t) y(t) \left[ \int_{h_2(t)}^t \frac{K_2+\alpha_2 x(s)}{1+x(s)}~d_s R_2(t,s) - y(t) \right]
\end{array}
\end{equation}

The following result generalizes \cite[Theorem 3.3.4, p. 193]{Gopalsamy}.

\begin{guess}
\label{theorem_appl2}
Suppose (a3)-(a6) hold and $\alpha_i>K_i$,  $i=1,2$.

Then there exists a unique positive equilibrium of \eqref{summary4}, and
all solutions of \eqref{summary4},\eqref{2star} converge to this equilibrium.
\end{guess}

A natural generalization of the results of the present paper would be to 
$n$-dimensional cooperative systems, as well as models with general nonlinear non-delay 
mortality
$$
\begin{array}{ll}
\displaystyle \frac{dx}{dt} & \displaystyle  = r_1(t) \left[ \int_{h_1(t)}^t f_1(y(s))~d_s R_1 (t,s) - g_1(x(t)) \right] 
 \vspace{2mm} \\
\displaystyle \frac{dy}{dt} & \displaystyle  = r_2(t) \left[ \int_{h_2(t)}^t f_2(x(s))~d_s R_2 (t,s) - g_2(y(t)) \right]
\end{array}
$$
In the one-dimensional case and monotone increasing $f_i$, such models have the same properties as equations 
with linear mortality functions $g_i(x)=b_ix$ \cite{Nonlin2013}.
So far we considered the case of the unique coexistence equilibrium; however, it would be interesting
to study multiple coexistence equilibria. For a single equation this investigation was implemented in  \cite{JMAA2014}.

\ack
The second author was partially supported by the NSERC grant RGPIN/261351-2010.
The authors are very grateful to the anonymous referee whose thoughtful suggestions significantly 
contributed to the present form of the paper.

\end{document}